\documentclass[11pt]{article}
\usepackage{amsthm}
\usepackage{youngtab}
\usepackage{xypic}
\usepackage{pstricks,epsfig}
\usepackage{pst-plot}
\usepackage{graphicx}
\usepackage[T1]{fontenc}
\usepackage[latin1]{inputenc}
\usepackage{amsfonts,amssymb,amsxtra}

\usepackage{latexsym}
 
\newcommand{\el}{\par \mbox{} \par \vspace{-0.5\baselineskip}}

\newcounter{amoi}
\newtheorem{theo}{Theorem}[section]
\newtheorem{leme}[theo]{Lemma}
\newtheorem{prop}[theo]{Proposition}

\newtheorem{defi}[theo]{Definition}

\newtheorem{fait}[theo]{Fact}
\newtheorem{exem}[theo]{Example}
\newcommand{\noi}{\noindent}
\newtheorem{rem}[theo]{Remark}

\newenvironment{preu}{{\flushleft \bf Proof :}}{\hfill $\square$ \vspace{2mm}}

\def\essai {\cdot\!\ \!\ \!\!\!\cdot\!\!\ \!\ \!\!\cdot}

\def \D {{\underline{D}}}
\def \im {{\mathrm{Im}}}
\def \p {{\mathbb{P}}}
\def \ker {{\mathrm{Ker}}}
\def\N{{\mathbb{N}}}

\def \Z {{\mathbb{Z}}}
\def\boxit#1#2{\setbox1=\hbox{\kern#1{#2}\kern#1}
\dimen1=\ht1 \advance\dimen1 by #1 \dimen2=\dp1 \advance\dimen2 by #1
\setbox1=\hbox{\vrule height\dimen1 depth\dimen2\box1\vrule}
\setbox1=\vbox{\hrule\box1\hrule}
\advance\dimen1 by .4pt \ht1=\dimen1
\advance\dimen2 by .4pt \dp1=\dimen2 \box1\relax}

\def\adots{\mathinner{\mkern2mu\raise1pt\hbox{.}
\mkern3mu\raise4pt\hbox{.}\mkern1mu\raise7pt\hbox{.}}}
\def\<{\langle\,}
\def\>{\,\rangle}

\def\codim{{\rm Codim}}

\def \cN {{{\mathcal{N}}}}
\def \cF {{{\mathcal{F}}}}
\def \tcN {{{\widetilde{\mathcal{N}}}}}

\def \osp {\mathfrak{osp}}

\begin{document}

\title{Study of  some orthosymplectic Springer fibers}
\author{S. Leidwanger \& N. Perrin}

\maketitle

\begin{abstract} 
We decompose the fibers of the Springer resolution for the odd nilcone of the 
Lie superalgebra $\osp(2n+1,2n)$ into locally closed subsets. We use this 
decomposition to prove that almost all fibers are connected. However, in 
contrast with the classical Springer fibers, we prove that the fibers can be
disconnected and non equidimensional.
\end{abstract}

\section*{Introduction}
As for classical Lie algebras, the odd nilpotent cone $\cN_1$ of the
Lie superalgebra $\osp(2n+1,2n)$ has a natural resolution (cf. \cite{GL} and 
section \ref{rappels}). We call it the \emph{Springer resolution} and denote 
it by $\pi:\tcN_1\to\cN_1.$
The purpose of the present paper is to describe some properties of the
fibers of $\pi$. On the one hand, some results, true for Lie algebras, are no 
longer true in the Lie superalgebra setting. Indeed, the study of explicit 
examples leads to the following proposition.

\begin{prop}
The fibers of $\pi$ are, in general,  neither connected nor
equidimensional. In particular, the variety $\cN_1$  is not normal.
\end{prop}

On the other hand we give for the orthosymplectic Lie superalgebra, as in the 
Lie algebra setting, a decomposition of the fiber of $\pi$ into locally closed 
subsets (see Theorem
\ref{theo-decoup}). These subsets do not have the same dimension in
general and their closures are  not always  irreducible components
of the Springer fiber. We use this decomposition  to prove the
following result. 

\begin{theo}
\label{main}
Let $\mathcal{O}_1$ be the unique codimension 1 orbit in $\cN_1$ and
let $X\in\cN_1$. The fiber $\pi^{-1}(X)$ is connected if and only if
$X\not\in\mathcal{O}_1$.
\end{theo}

As an application of our decomposition, we also describe explicitely
the fibers of $\pi$ above the three bigger non dense orbits of $\cN_1$
under the action of the orthosymplectic group. For the unique codimension 3
orbit, the fiber is not equidimensional and has irreducible components
of dimension 1 and 2 (see Proposition \ref{prop-3}).
   
\section{Decomposition of the fiber}
\label{rappels}

\subsection{Some basic facts about odd nilpotent orbits}
Let $V=V_0\oplus V_1$ be a  $\Z/2\Z$-graded vector space of
super-dimension $(2n+1,2n)$ and  equipped with a bilinear
super-symmetric form  $B$. This means that the spaces  $V_0$ and
$V_1$ are orthogonal, of respective dimensions  $2n+1$ and $2n$  and
the restrictions  $\varphi:V_0\to V_0^\vee$ and $\psi:V_1\to V_1^\vee$ of
$B$ to $V_0$ and $V_1$ are non degenerate respectively symmetric and alternate.
Let $\osp(2n+1,2n)$ be the Lie superalgebra, called the orthosymplectic
Lie superalgebra, consisting in endomorphisms of $V$ which preserve the
bilinear super-symmetric form and the graduation (for more details,
see for example \cite{GL}).

The group ${\rm O}(V_0,\varphi)\times {\rm Sp}(V_1,\psi)$ acts on
$\osp(2n+1,2n)$. It is called the orthosymplectic group and we denote it by 
$G_0\times G_1$. 
 Let $u$ be an element of  ${\rm Hom}(V_0,V_1)$,
 we define $u^*\in {\rm Hom}(V_1,V_0)$ by 
$u^*=\varphi^{-1}\circ u^t\circ \psi$.
An endomorphism $X=(u,u^*)$ of $V$,  with $u\in{\rm Hom}(V_0,V_1)$  is
called odd orthosymplectic and belongs to $\osp(2n+1,2n)$. Any degree
1 element in $\osp(2n+1,2n)$ preserving the $\Z/2\Z$-graduation is of this shape.

The set of odd nilpotent orthosymplectic endomorphisms of ${\rm
  End}(V)$ is a cone denoted by $\cN_1$. A $G_0\times G_1$-equivariant
resolution of the singularities of $\cN_1$ is constructed in
\cite{GL}. Let us describe this resolution.

Let  $B_0$ and $B_1$ be Borel subgroups of $G_0$ and $G_1$. For
$X=(u,u^*)\in \cN_1$ we denote by ${B}_u={B}_X$ the set of  pairs of
isotropic complete flags $((E_i)_{i\in[1,n]},(F_j)_{j\in[1,n]})$ in $G_0/B_0\times G_1/B_1$,
with $E_i\subset V_0$,  $\dim E_i=i$, $F_j\subset V_1$, $\dim F_j=j$ and for all $i\in[1,n]$
\begin{equation}
  \label{eq-fibre}
X(E_i)=u(E_i)\subset F_{i-1} \textrm{ and } X(F_i)=u^*(F_i)\subset E_i. 
\tag{\dag}
\end{equation}
By analogy with the classical case,   the variety ${B}_u={B}_X$
corresponding to $X=(u,u^*)\in\cN_1$ is called Springer fiber since it
can be identified with the fiber above $X$ of the resolution $\tcN_1$  of  the
singularities of $\cN_1$: 
$$\tcN_1=\{(u,(E_i)_{i\in[1,n ]},(F_j)_{j\in[1,n]})\in {\rm Hom}(V_0,V_1)
\times G_0/B_0\times G_1/B_1\ /\ ((E_i),(F_j))
\in B_u\}.$$
The map $\pi:\tcN_1\to \cN_1$ is the first projection, the second one
$p_2$ realizes $\tcN_1$ as a vector bundle above $G_0/B_0\times
G_1/B_1$.  

\vskip 0.5 cm

H.P. Kraft and  C. Procesi \cite{KP} proved that the odd
nilpotent orthosymplectic orbits under $G_0\times G_1$-action of $\osp(m,2n)$ (where $m\in \{2n+2, 2n+1, 2n, 2n-1, 2n-2\}$)  are parametrized  by marked Young diagrams of size $m+n$
(see Fulton \cite{Ful97} for more details on Young diagrams). We
recall their results and specify that the diagrams are  written in the
french way. 

\begin{defi}
\label{def1}
\label{indecom} 
(\i) A \emph{marked diagram} of  size $(m,n)$ is a Young diagram of size $m+n$
in which there are $m$ boxes labelled with $0$ and $n$ boxes labelled with $1$. The labels in
the same line alternate. 

(\i\i) A line beginning with  $\epsilon\in\{0,1\}$ is said to be of
parity $\epsilon$. 

(\i\i\i) A marked diagram $D$ is  called  \emph{indecomposable} if it
has one of the following shapes:

\hskip 1 cm 1. an even line of length  $4p+1$,

\hskip 1 cm 2. an odd line of length $4p-1$,

\hskip 1 cm 3. two even lines of length $4p-1$,

\hskip 1 cm 4. two odd lines of length $4p+1$,

\hskip 1 cm 5. two lines, one even, one odd of length $2p$.

(\i v) A marked  diagram is  \emph{admissible} if it is the union of 
indecomposable diagrams.
\end{defi}

\begin{prop} \cite{KP} There is a bijective correspondence between odd 
nilpotent orbits of $\osp(m,2n)$ and admissible diagrams of size $(m,n)$.
\end{prop}

An easy consequence of the above correspondence is the following fact.

\begin{fait}
\label{fait-noyau}
\label{fait-noyau-gen}
 Let $X\in\cN_1$, $D$ be its associated diagram by the previous correspondence. The dimension  of the space $Ker X\cap \im
 X^{k-1}\cap V_\epsilon$ is the number of lines of parity $\epsilon$ and of
 length at least $k$. The super form is trivial on this space except for
 $k=1$. In this case  its rank is the number of marked lines of parity 
$\epsilon$ and length 1.  
\end{fait}

\subsection{Slicing  of diagrams}
 An admissible diagram $D$ is said to have a parity  if its size is $(2n +
(-1)^\epsilon, 2n)$ for some $n$ and  $\epsilon \in \{0,1\}$ ; the parity of $D$ is then 
$\epsilon$.
(If $D$ is a line of the given size this coincides with the previous notion of parity.)

\begin{defi}
  Let $D$ be an admissible diagram having a parity. We call  \emph{admissible subdiagram
    of $D$} any subdiagram $D'$ of $D$ such that

\hskip 0.3 cm 1. $D'$ is admissible, 

\hskip 0.3 cm 2. $D'$ has two boxes less than $D$,

\hskip 0.3 cm 3. $D'$ and $D$ have different parities,

\hskip 0.3 cm 4. the boxes of $D\setminus D'$ are at the beginning and
at the end  of lines of same length. 
\end{defi}

\begin{leme}
\label{lemm-deux-cas}
Let $D$ be an admissible diagram of parity $\epsilon$ and $k$ be an integer 
such that $D$ has at least one line of length $k$ and parity $\epsilon$. Then 
there exist at least one and at most two admissible subdiagrams $D'$ of 
$D$ such that the boxes in $D\setminus D'$ lie on lines of length $k$. 
  \end{leme}

\begin{preu}
The two boxes are removed either on the same line or on two
different lines. The different cases  are as follows.

$\bullet$ 
If $k$ is even, the admissibility of $D$ and $D'$ (definition
    \ref{indecom}-(5)) implies  that  the  number of lines of length
    $k$ of $D$ and $D'$ is even and then two lines of length $k$
    change size i.e. the   set $D\setminus D'$ is on two different
    lines. 

$\bullet$ If $k=4p+(-1)^\epsilon$,  the admissibility of $D$ and $D'$
  (definition \ref{indecom}-(1)-(2)) implies that we can have any
  number of lines of length $k$ in $D$ 
  and of lines of length $k-2$ in $D'$. We can choose $D\setminus D'$
on the same line or on  two different lines if there exist two lines of 
length $k$.

$\bullet$ If $k=4p-(-1)^\epsilon$, the
  admissibility of $D$ and $D'$ (definition \ref{indecom}-(3)-(4))
  implies that the number of lines of length $k$ in $D$ has to  be even
  and then that two lines of length $k$ change size i.e. the set
  $D\setminus D'$ is on two different  lines.  
\end{preu}

\begin{defi}
  Let $D$ be an admissible diagram of size $(2n+1,2n)$ and $D_0=D$.
A sequence $(D_i)_{i\in[0,2n]}$ such that $D_{i}$ is a admissible subdiagram 
of $D_{i-1}$ for $i\in[1,2n]$, is called an
\emph{admissible slicing} of $D$ and is denoted by 
$\underline{D}$. We denote by $\mathcal{A}(D)$ the set of admissible
slicings of $D$.  
\end{defi}

\begin{exem} We give below the unique two admissible slicings of a diagram.
$$
\begin{array}{l}
  \begin{array}[u]{lllllllll}
  \young(0,101,01010)&\textrm{\raisebox{3ex}{$\supset$}}&
\young(0,101,101)&\textrm{\raisebox{3ex}{$\supset$}}&
  \young(0,0,101)& \textrm{\raisebox{3ex}{$\supset$}} &
\textrm{\raisebox{2.5ex}{$\young(101)$}}& \textrm{\raisebox{3ex}{$\supset$}}&
\textrm{\raisebox{2.5ex}{$\young(0)$}} \\
\end{array}\\
\\
\begin{array}[u]{lllllllll}
  \young(0,101,01010)&\textrm{\raisebox{3ex}{$\supset$}}&
  \young(0,101,101)&\textrm{\raisebox{3ex}{$\supset$}}&
  \young(0,10,01)& \textrm{\raisebox{3ex}{$\supset$}} &
\young(0,1,1)& \textrm{\raisebox{3ex}{$\supset$}}&
\textrm{\raisebox{2.5ex}{$\young(0)$ .}} \\
\end{array}\\
\end{array}$$
\end{exem}

\subsection{Locally closed subsets in the fiber}

Let $X$ be a nilpotent element in $\osp(2n+(-1)^\epsilon,2n)$. We denote by  
$\mathcal{K}_\epsilon(X,k)$ the set of isotropic points of
$\p(V_\epsilon\cap \ker X\cap\im X^{k-1}) \setminus\p(V_\epsilon\cap
\ker X\cap\im X^{k})$ and by $\ell_\epsilon(X,k)$ the number of lines
of length $k$ and parity $\epsilon$ in the diagram associated to $X$.

\begin{rem}
  \label{rem-dim-K}
(\i) When $\epsilon=1$ or when $\epsilon=0$ and $\ell_0(X,1)=0$, all
  the points of the projective space $\p(V_\epsilon\cap \ker
  X\cap\im X^{k-1})$ are isotropic.

(\i\i) The variety $\mathcal{K}_0(X,1)$ is empty when $\ell_0(X,1)=1$.
For the other cases we have
$$\dim\mathcal{K}_\epsilon(X,k)=\left\{
  \begin{array}{ll}
    \sum_{i\geq k}\ell_\epsilon(X,i)-2& \textrm{if $\epsilon=0$,
      $k=1$ and $\ell_0(X,1)>1$}\\
    \sum_{i\geq k}\ell_\epsilon(X,i)-1& \textrm{otherwise.}\\
  \end{array}\right.$$

(\i\i\i)The variety $\mathcal{K}_\epsilon(X,k)$ is irreducible (therefore 
connected) except when $\epsilon=0$, $k=1$ and $l_0(X,1)=2$. In
  this last case, the isotropic locus is the union of two 
hyperplanes $H$ and $H'$ in $\p(\ker X)$ and we have $\mathcal{K}_0(X,1)=(H\cup
H')\setminus(H\cap H')$.
\end{rem}

\begin{exem}
Let $X$ be a nilpotent element in $\osp(5,4)$ with associated diagram
as follows:
$$\textrm{\raisebox{3ex}{$D=$ }}\young(01,10,01010)\ .$$
Its kernel has dimension 3 and its intersection with $V_0$ is totally
isotropic of dimension 2. We have then $\dim\mathcal{K}_0(X,2)=1$ and
$\dim\mathcal{K}_0(X,5)=0$.
\end{exem}

\begin{prop}
\label{prop-cle}
 Let $X\in\osp(2n+(-1)^\epsilon,2n)$ be a nilpotent element with
 associated diagram $D$, let $x\in\mathcal{K}_\epsilon(X,k)$ and let
 $y$ be such that $X^{k-1}(y)=x$. The restriction $X\vert_{x^\perp/x}$
 of $X$ to $x^\perp/x$ lies in $\osp(2n-1,2n-2\epsilon)$ and its
 orbit under the corresponding orthosymplectic group is associated to
 the admissible diagram 
\begin{itemize}
\item obtained by removing two boxes in a line  of length $k$ of $D$
  if  $B(x,y)\neq0$,

\item obtained by removing two boxes in two different lines of length  $k$ of
  $D$ if $B(x,y)=0$.
 \end{itemize}
\end{prop}

\begin{preu}
We first notice that if $X$ is an odd orthosymplectic nilpotent
element, its diagram is entirely determined by the dimensions of $\ker
X^a\cap V_0$  and  $\ker X^a\cap V_1$ for all $a$. We therefore
have to compute $\dim (\ker(X\vert_{x^\perp/x}^a))$.

We first determine $\ker(X^a)\cap x^\perp$. If $z\in\ker X^a$ with $a\leq k-1$, 
we have the equality
$B(x,z)=B(X^{k-1}(y),z)=B(y,X^{k-1}(z))=0$ i.e. 
$\ker X^a\subset x^\perp$ for $a\leq k-1$. But 
$\ker X^{k}\not\subset x^\perp$ otherwise we get $(\im
X^{k})^\perp=\ker X^{k} \subset x^\perp$ and  $x\in\im X^{k}$, a
contradiction with $x\in\mathcal{K}_\epsilon(X,k)$. Recall that $x\in
V_\epsilon$  we thus obtain the equalities:
$$\dim(\ker X^a\cap x^\perp\cap V_\epsilon)=\left\{
  \begin{array}{ll}
    \dim(\ker X^a\cap V_\epsilon)& \textrm{if $a\leq k-1$}\\
    \dim(\ker X^a\cap V_\epsilon)-1& \textrm{if $a>k-1$,}\\
  \end{array}\right.$$
$$\dim(\ker X^a\cap x^\perp\cap V_{1-\epsilon})=\dim(\ker X^a\cap 
V_{1-\epsilon})\ \textrm{for all $a\in\N$.}$$

Let us consider $Y=X\vert_{x^\perp}$ which is possible since 
$\im X=\ker X^\perp\subset x^\perp$. We compute the dimensions of 
$\ker Y^a=\ker X^a\cap x^\perp$. Set
$Z=X\vert_{x^\perp/x}$. We compute $\dim\ker Z^a$ using $\dim\ker Y^a$. 
By definition, we have $\ker Z^a=(Y^a)^{-1}(\langle
x\rangle)/x$, we obtain
$$\dim((Y^a)^{-1}(\langle x\rangle)=\left\{
  \begin{array}{ll}
    \dim\ker X^a+1& \textrm{si $a\leq k-1$}\\
    \dim\ker X^a& \textrm{si $a>k-1$.}\\
  \end{array}\right.$$
If $a\leq k-1$, the element $X^{k-1-a}(y)$ is in $(Y^a)^{-1}(\langle
x\rangle)$ but not in $\ker X^a$. Therefore we have the equalities:
$$\dim((Y^a)^{-1}(\langle x\rangle)\cap x^\perp)=\left\{
  \begin{array}{ll}
    \dim\ker X^a+1& \textrm{if $a<k-1$,}\\
    \dim\ker X^a+1& \textrm{if $a=k-1$ and $B(x,y)=0$,}\\    
\dim\ker X^a& \textrm{if $a=k-1$ and $B(x,y)\neq0$,}\\
    \dim\ker X^a& \textrm{si $a>k-1$.}\\
  \end{array}\right.$$
We have $\dim \ker Z^a=\dim((Y^a)^{-1}(\langle x\rangle)\cap
x^\perp)-1$ and for $x\in V_\epsilon$, we have the inclusion
$(Y^a)^{-1}(\langle x\rangle)\subset V_{\epsilon+a\mod 2}$.
We obtain then the following dimensions.

$$\dim(\ker Z^a\cap V_\epsilon)=\left\{
  \begin{array}{ll} 
    \dim(\ker X^a\cap V_\epsilon)-1& \textrm{if $a<k-1$ and $a$ is odd,}\\
    \dim(\ker X^a\cap V_\epsilon)& \textrm{if $a<k-1$ and $a$ is even,}\\
    \dim(\ker X^a\cap V_\epsilon)-1& \textrm{sif $a=k-1$, $k-1$ is odd 
and $B(x,y)=0$,}\\    
    \dim(\ker X^a\cap V_\epsilon)& \textrm{if $a=k-1$, $k-1$ is even and 
$B(x,y)=0$,}\\    
\dim(\ker X^a\cap V_\epsilon)-1& \textrm{if $a=k-1$ and $B(x,y)\neq0$,}\\
    \dim(\ker X^a\cap V_\epsilon)-1& \textrm{if $a>k-1$.}\\
  \end{array}\right.$$
$$\dim(\ker Z^a\cap V_{1-\epsilon})=\left\{
  \begin{array}{ll}
    \dim(\ker X^a\cap V_{1-\epsilon})+1& \textrm{if $a<k-1$ and $a$ is 
odd,}\\
    \dim(\ker X^a\cap V_{1-\epsilon})& \textrm{if $a<k-1$ and $a$ is 
even,}\\
    \dim(\ker X^a\cap V_{1-\epsilon})+1& \textrm{if $a=k-1$, $k-1$ is odd and  
$B(x,y)=0$,}\\    
    \dim(\ker X^a\cap V_{1-\epsilon})& \textrm{if $a=k-1$, $k-1$ is even
and $B(x,y)=0$,}\\    
\dim(\ker X^a\cap V_{1-\epsilon})& \textrm{if $a=k-1$ and $B(x,y)\neq0$,}\\
    \dim(\ker X^a\cap V_{1-\epsilon})& \textrm{if $a>k-1$.}\\
  \end{array}\right.$$
The diagrams determined by these  dimensions are the required diagrams.
\end{preu}

Let $X$ be an odd nilpotent element in $\osp(2n+1,2n)$ and let $D(X)$
be its associated diagram. Let $\D=(D_i)_{i\in[0,2n]}$ be an
admissible slicing of $D(X)$. We denote by $X_{2i-1}$ (resp. $X_{2i}$)
the restriction of  $X$ to
${(E_i^\perp/E_i)\cap(F_{i-1}^\perp/F_{i-1})}$ (resp. to
${(E_i^\perp/E_i)\cap(F_{i}^\perp/F_{i})}$). 

\begin{defi}
We define the subset $B_X(\D)$ of the fiber $(X,B_X)$ by
$$B_X(\D)=\{((E_i)_{i\in[1,n]},(F_j)_{j\in[1,n]})\in p_2(\pi^{-1}(X))\ /\
D(X_{2i-1})=D_{2i-1}\ {\rm
    and}\ D(X_{2i})=D_{2i}\}.$$
\end{defi}

\begin{theo}
  \label{theo-decoup}
 Let  $X$ be an odd nilpotent element in $\osp(2n+1,2n)$. 
The subsets $B_X(\D)$ are locally closed in $B_X$ and we have
$$B_X=\coprod_{\D \in \mathcal{A}(D(X))}B_X(\D).$$ 
\end{theo}

\begin{preu}
Let  $X$ be an odd nilpotent element in $\osp(2n+1,2n)$,
we proceed by induction on the size of the diagram $D=D(X)$.
Denote by $\cF_k$ and $\cF'_k$ the varieties
$$
\begin{array}{l}
  \cF_k=\{((E_i)_{i\in[1,k]},(F_j)_{j\in[1,k]})\ /\ \textrm{  $E_i$ and 
$F_j$ satisfy the equation (\ref{eq-fibre})}\}  \textrm{ and } \\
\cF'_k=\{((E_i)_{i\in[1,k]},(F_j)_{j\in[1,k-1]})\ /\ \textrm{  $E_i$ and 
$F_j$ satisfy the equation (\ref{eq-fibre})}\}.\\
\end{array}$$ 
We have a sequence of morphisms
$\cF_n\to\cF'_n\to\cF_{n-1}\to\cdots\cF_1\to\cF'_1\to
\cF_0=\{{\rm pt}\}.$
The fiber of the morphism $\cF'_i\to\cF_{i-1}$
(resp. $\cF_i\to\cF'_i$) is given by the isotropic elements of
$\p(\ker Y)$ where $Y$  is the restriction of $X$ 
to ${(E_{i-1}\oplus F_{i-1})^\perp/(E_{i-1}\oplus F_{i-1})}$ (resp. to
${(E_{i}\oplus F_{i-1})^\perp/(E_{i}\oplus F_{i-1})}$). Those $Y$ are
orthosymplectic and their associated diagrams are as in  Proposition
\ref{prop-cle} (i.e. of size less than the size of $D$).
These fibrations are locally trivial.

If $D_1$ is obtained from $D$ by removing boxes on lines of length $k$, then
the fiber of the map $\cF'_1\to\cF_0=\{{\rm pt}\}$ is the locally closed
subset  $\mathcal{K}_0(X,k)$. 
 We then consider $X\vert_{E_1^\perp/E_1}$ and apply the induction hypothesis.
\end{preu}

\begin{rem}\label{rem-dim-fibre}
This result may remind the reader of results of Spaltenstein \cite{spal} and Van Leuwen 
\cite{vanL} for classical (types $A$, $B$, $C$ and $D$) Lie algebras. However 
in the Lie algebra setting the dimensions of the locally closed subsets 
obtained by admissible slicing are constant. The closure of these locally 
closed subsets are therefore the irreducible components of the Springer 
fiber. This does not happen in our situation.
\end{rem}

\begin{exem}
Let $X$ be a nilpotent element in $\osp(5,4)$ with associated diagram $D(X)$ 
as follows. The admissible slicings of  $D(X)$ are:
$$
\begin{array}{l}
\textrm{\raisebox{3ex}{$D(X)=$ }}
\young(01,10,01010)
\textrm{\raisebox{3ex}{  $\supset$ }}
\young(1,1,01010)
\textrm{\raisebox{3ex}{  $\supset$ }}
\textrm{\raisebox{2.6ex}{  $\young(01010)$ }}
\textrm{\raisebox{3ex}{  $\supset$ }}
\textrm{\raisebox{2.6ex}{  $\young(101)$ }}
\textrm{\raisebox{3ex}{  $\supset$ }}
\textrm{\raisebox{2.6ex}{  $\young(0)$ ,}}\\
\\
\textrm{\raisebox{3ex}{$D(X)=$ }}
\young(01,10,01010)
\textrm{\raisebox{3ex}{  $\supset$ }}
\young(01,10,101)
\textrm{\raisebox{3ex}{  $\supset$ }}
\young(0,0,101)
\textrm{\raisebox{3ex}{  $\supset$ }}
\textrm{\raisebox{2.6ex}{  $\young(101)$ }}
\textrm{\raisebox{3ex}{  $\supset$ }}
\textrm{\raisebox{2.6ex}{  $\young(0)$ ,}}\\
\\
\textrm{\raisebox{3ex}{$D(X)=$ }}
\young(01,10,01010)
\textrm{\raisebox{3ex}{  $\supset$ }}
\young(01,10,101)
\textrm{\raisebox{3ex}{  $\supset$ }}
\young(0,01,10)
\textrm{\raisebox{3ex}{  $\supset$ }}
\young(0,1,1)
\textrm{\raisebox{3ex}{  $\supset$ }}
\textrm{\raisebox{2.6ex}{  $\young(0)$ .}}
\end{array}
$$
The fiber $(X,B_X)$ is then an union of three components
$(X,B_X(\overline{D}))$ of respective  dimensions $2$, $1$, $1$. 
To determine the irreducible components of the fiber, we have to find the
locally closed subsets of dimension $1$  belonging to the closure of
the locally closed subset of dimension $2$. 
\end{exem}

\begin{exem}
Let $X\in \cN_1$ such that $D(X)$ is a hook. One easily checks that
the locally closed subsets of $\pi^{-1}(X)$ are equidimensional,
therefore $B_X$ is equidimensional with irreducible components indexed
by ${\cal A}(D(X))$. Let $p$ be an even integer, the dimension of
$B_X$ is equal to:

\begin{itemize}
\item $\frac{p^2}{2}$ if $D(X)$ has an even line of length $4n+1-2p$ and
  $2p$  lines of length 1 ($p$ even, $p$ odd), 
\item $\frac{p(p-2)}{2}+1$ if $D(X)$ has an even line of length
  $4n+1-2p$  and $2p-2$ lines of length 1 ($p$ odd, $p-2$ even), 
\item $\frac{p(p-2)}{2}$ if $D(X)$ has an odd line of length
  $4n+3-2p$ and  $2p-2$ lines of length 1 ($p$ even, $p-2$ odd),  
\item $\frac{p^2}{2}$ if $D(X)$ has an odd line of length $4n+1-2p$ and
  $2p$  lines of length 1 ($p$ even, $p$ odd).
\end{itemize}
\end{exem}

\section{Connectedness of the fibers?}

In this section we prove Theorem \ref{main} \emph{i.e.} we determine
for which element $X\in\cN_1$, the fiber $\pi^{-1}(X)$ is connected.
Recall that  $\mathcal{O}_1$ is the odd nilpotent subregular  orbit
\emph{i.e.} of codimension 1 in $\osp(2n+1,2n)$ (the fiber over this
orbit is disconnected, see Proposition \ref{prop-1}). 

\begin{theo}
 Let $X\in\cN_1$,  the fiber $\pi^{-1}(X)$ is connected if and only if
 $X\not\in\mathcal{O}_1$.
\end{theo}

\begin{preu}
 We proceed by induction on the size of the diagram associated to $X$
 and use the sequence of morphisms described in the proof of Theorem
 \ref{theo-decoup}:
$$\cF_n\to\cF'_n\to\cF_{n-1}\to\cdots\cF_1\to\cF'_1\to
\cF_0=\{{\rm pt}\}.$$
Let $K_0$ be the closure of $\mathcal{K}_0(X,1)$. The map
$p:\pi^{-1}(X)\to \cF'_1$ 
takes values in $K_0$  and is surjective. Let  $\D$ be an admissible
slicing of $ D$, then $p(B_X(\D))$ is locally closed in
$K_0$. If we consider every admissible slicing we obtain a
stratification of $K_0$. 

We prove that every point of a locally closed subset is connected, by
a curve, to a point of the special locally closed subset i.e. admiting
the most special admissible slicing (this slicing is obtained by
removing, when there is a choice, boxes on the longest lines). We then
prove that the fiber is connected if the corresponding orbit is not
$\mathcal{O}_1$. 

We first reduce to the case where $E_1\in K_0$ is in the smallest
stratum i.e. $E_1\subset\ker X\cap \im X^{k-1}$ for $k$ the maximal
length of an even line of $D$. By Remark \ref{rem-dim-K}, as soon as
$\dim K_0>0$, the closed subset $K_0$ is connected. Therefore, the 
surjectivity of $p$ on $K_0$ gives us a curve in $\pi^{-1}(X)$ such that
the image of the generic point is any point of  $K_0$ (for example
$E_1$) while the image of the special one is a point of $\p(\ker X\cap
\im X^{k-1})$. Remark also (see Remark \ref{rem-dim-K} again), that
$\dim K_0=0$ only for the orbit $\mathcal{O}_1$. 

Using this argument recursively, we are reduced to proving that if
$\D$ is the special slicing described above, the locally closed
subset $B_X(\D)$ is connected.

Let  $\D$ be the special admissible slicing of $D$.
To choose  $E_1$, we choose an isotropic element of 
$\p(\ker X\cap\im X^{k-1}\cap V_0)$ where $k$ is the maximal length of
an even line. If  in $D$ there is an even line of length at least  $2$,
then the quadratic form   $\varphi$ restricted to $\ker X\cap\im
X^{k-1}\cap V_0$ is trivial and we choose any point of 
 $\p(\ker X\cap\im X^{k-1}\cap V_0)$ which is connected. 
If not then we have only even lines of length  1. We cannot have only
one such line (otherwise we would be in $\osp(1,0)$!) and if there
are at least 3 of them, the set of isotropic elements of 
$\p(\ker X\cap\im X^{k-1}\cap V_0)$ is a quadric of dimension at least
1 therefore connected. It remains the case of exactly two even lines of
length 1, we notice that the diagram is then the diagram associated to
$\mathcal{O}_1$ 
$$\young(0,0,101\essai 101)\ ,$$
which is not allowed. Next we determine $F_1$. To do this we take an
element in a projective space (the isotropic condition is always
verified). 

We now have to verify  by induction that the diagram
corresponding to $\mathcal{O}_1$ can appear in the special admissible
slicing $\D$ of $D$ only if  $D$ itself corresponds to $\mathcal{O}_1$.  
If the diagram corresponding to  $\mathcal{O}_1$  is one of the $D_i$
of $\D$ for  $i>0$, then $D_{i-1}$ has one of the following two forms:
$$\young(0,0,1,1,101\essai 101)\
\textrm{\raisebox{6ex}{ or  }}\textrm{\raisebox{2.6ex}{
    $\young(01,10,101\essai 101)$ .}}$$
We notice that the length of the first line is odd, it has then to be
at least 3 which means that the choice of $D_i$ was not the
most special, a contradiction.
\end{preu} 

\section{Fibers for the orbits $\mathcal{O}_1$, $\mathcal{O}_2$ and 
$\mathcal{O}_3$}

We study in this section the fibers of $\pi$ above the non dense
orbits of maximal dimension $\mathcal{O}_i$  with $\codim\mathcal{O}_i=i
\in[1,3]$. The diagram of the orbit $\mathcal{O}_1$ is given in the previous 
proof.
\begin{prop}
  \label{prop-1}
For  $X\in\mathcal{O}_1$, the fiber  $\pi^{-1}(X)$ is the disjoint
union of two points.
\end{prop}

\begin{preu}
The decomposition in locally closed subsets has only one element.
Moreover, this locally closed subset is of dimension $0$ and has two
connected components. 
\end{preu}

The diagrams of $\mathcal{O}_2$ and $\mathcal{O}_3$ are  as follows
(the first lines have length  $4n-3$):  
$$\young(0,101,0101\essai 1010)\qquad \textrm{\raisebox{3ex}{ and
  }}
\qquad \young(10,01,0101\essai 1010).$$

\begin{prop}
  \label{prop-2} Let $X$ be in  the orbit
  $\mathcal{O}_2$ and assume that $n\geq2$.

(\i) The fiber  $\pi^{-1}(X)$ is non-reduced everywhere.

(\i\i) The reduced fiber, denoted by $\pi^{-1}(X)_{\textrm{r{e}d}}$,
  is the union of $2n-1$ irreducible components
  $(C_i)_{i\in[1,2n-1]}$, all isomorphic to $\p^1$, such that
\begin{itemize}
\item the components  $(C_{2i-1})_{i\in[1,n-2]}$ (resp. 
$(C_{2i})_{i\in[1,n-2]}$) form a chain that meets transversally the  
component $C_{2n-3}$ in $x$ (resp. in $y$ distinct from $x$),
\item the components $C_{2n-2}$ and  $C_{2n-1}$ meet transversally the
  component $C_{2n-3}$ in two distinct points (also distinct from
  $x$ and $y$).
\end{itemize}
The dual graph of  $\pi^{-1}(X)_{\textrm{r{e}d}}$ is then the
following (the left branches have length $n-2$).\\ 
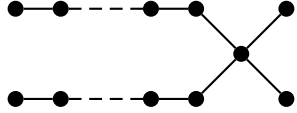
\centerline{\begin{pspicture*}(0,1)(6.0,3)

\psellipse[fillstyle=solid,fillcolor=black](1.2,2.4)(0.106,0.106)
\psellipse[fillstyle=solid,fillcolor=black](1.2,1.2)(0.106,0.106)
\psellipse[fillstyle=solid,fillcolor=black](1.8,2.4)(0.106,0.106)
\psellipse[fillstyle=solid,fillcolor=black](1.8,1.2)(0.106,0.106)
\psellipse[fillstyle=solid,fillcolor=black](3.0,2.4)(0.106,0.106)
\psellipse[fillstyle=solid,fillcolor=black](3.0,1.2)(0.106,0.106)
\psellipse[fillstyle=solid,fillcolor=black](3.6,2.4)(0.106,0.106)
\psellipse[fillstyle=solid,fillcolor=black](3.6,1.2)(0.106,0.106)
\psellipse[fillstyle=solid,fillcolor=black](4.2,1.8)(0.106,0.106)
\psellipse[fillstyle=solid,fillcolor=black](4.8,2.4)(0.106,0.106)
\psellipse[fillstyle=solid,fillcolor=black](4.8,1.2)(0.106,0.106)
\psline(1.306,2.4)(1.694,2.4)
\psline(3.106,2.4)(3.494,2.4)
\psline(3.675,2.325)(4.125,1.875)
\psline(4.275,1.875)(4.725,2.325)
\psline(4.725,1.275)(4.275,1.725)
\psline(4.125,1.725)(3.675,1.275)
\psline(3.494,1.2)(3.106,1.2)
\psline(1.694,1.2)(1.306,1.2)
\psline[linestyle=dashed](1.906,2.4)(2.894,2.4)
\psline[linestyle=dashed](1.906,1.2)(2.894,1.2)

\end{pspicture*}
}
\end{prop}

\begin{preu}
(\i) To find the one dimensional subspaces in  $\ker X\cap V_0$, we
  have to look for the isotropic points of $\p(\ker X\cap
  V_0)$. According to Fact \ref{fait-noyau-gen}, the vector space
  $\ker X \cap V_0$ 
  has dimension $2$ and the quadratic form has rank one, the
  unique solution is then a double point and the fiber is not
  reduced. 

(\i\i) Let us consider the decomposition into locally closed subsets,
  obtained in Theorem \ref{theo-decoup}. Let $\D$ be an admissible
  slicing of  $D$, the diagram associated to $\mathcal{O}_2$. We have
  $D_0=D$ 
and  $D_1$ has three lines: two odd, one of length $4n-5$, the other of
length $ 3$ and one even of length $1$.  

There are two cases for $D_2$. Let us denote by  $D_2^g$ the
general one and by $D_2^s$ the special one. $D_2^g$ has  three lines:
one odd of  length $4n-5$, two even of length $1$. The diagram  
$D_2^s$ has three lines: one even of length  $4n-7$, another odd of
length $3$ and  the last even of length one. 
$$\textrm{\raisebox{3ex}{$D_1=$ }}\young(0,101,101\essai 101)\qquad 
\textrm{\raisebox{3ex}{$D_2^g=$ }}\young(0,0,101\essai 101)\qquad 
\textrm{\raisebox{3ex}{$D_2^s=$ }}\young(0,101,01\essai 10).$$
In the case $D_2^g$, we have an unique choice for  $D_3$ (one line of
length $4n-5$).  The diagrams $D_i$ for $i\geq3$ are then fixed (they
have alternatively one odd line or one even line).
In the case $D_2^s$, we recognize the diagram associated to the orbit
$\mathcal{O}_2$ in $\osp(2n-1,2n-2)$. An easy induction on
$n$ proves that there are $n$ admissible slicings.

We proceed by induction on $n$ to prove the proposition.
 To start the induction, we study $\osp(7,6)$. From the above, there
 are three admissible slicings. We denote them by $\D$, $\D'$ 
and $\D''$. Let us now describe explicitly the associated locally closed 
subsets.  To do this, we fix a representative $(u, u^*)$ of $\mathcal{O}_2$ 
such that the matrix of  
$u\in\textrm{Hom}(V_0,V_1)$ is \\
\begin{minipage}[b]{0.5\linewidth}
\centering
$$\left(\begin{array}{ccccccc}
  0&1&0&0&0&0&0\\
  0&0&1&0&1&0&0\\
  0&0&0&1&0&0&0\\
  0&0&0&0&0&0&0\\
  0&0&0&0&0&1&0\\
  0&0&0&0&0&0&1\\
\end{array}\right)$$
\end{minipage}
\hspace{-1cm}
\begin{minipage}[b]{0.5\linewidth}
in the bases $(e_i)_{i\in[1,7]}$
and $(f_i)_{i\in[1,6]}$ of $V_0$ and  $V_1$  such that 
$(e_i,e_j)=\delta_{i,8-j}$ and $(f_i,f_j)=\delta_{i,7-j}$. 
\label{fig:figure2}
\end{minipage}

The locally closed subsets are
$$\begin{array}{l}
  B_X(\D)=\left\{
    \begin{array}{ll}
\langle e_1\rangle
\subset
\langle e_1,\alpha e_2+\beta e_\pm\rangle \subset
\langle e_1,e_2,e_\pm\rangle\\
\langle \alpha f_1+\beta f_3\rangle\subset
\langle f_1,f_3\rangle\subset\langle f_1,f_2,f_3\rangle\end{array},
\textrm{ $[\alpha:\beta]\in\p^1$},\ \beta\neq 0
\right\}\\
B_X(\D')=\left\{
\begin{array}{ll}
  \langle e_1\rangle
\subset
\langle e_1,e_2\rangle\subset \langle e_1,e_2,e_\pm(\gamma,\delta)\rangle\\
\langle f_1\rangle\subset
\langle f_1, \gamma f_2+\delta f_3\rangle\subset
\langle f_1,f_2,f_3\rangle
\end{array},
\textrm{ $[\gamma:\delta]\in\p^1$},\
\displaystyle{\frac{\gamma^2}{2}+\delta^2\neq0}
\right\}\\ 
B_X(\D'')=\left\{
\begin{array}{ll}
  \langle e_1\rangle \subset 
\langle e_1,e_2 \rangle\subset
\langle e_1,e_2,e_\pm(\gamma,\delta)\rangle\\
\langle f_1\rangle\subset
\langle f_1, \gamma f_2+\delta f_3\rangle\subset 
\langle f_1, \gamma f_2+\delta f_3, f_{\gamma,\delta}(\zeta,\eta)
\rangle
\end{array},
\displaystyle{\frac{\gamma^2}{2}+\delta^2=0},\ [\zeta:\eta]\in\p^1\right\}\\ 
\end{array}$$
where
$e_\pm=e_4\pm\frac{\sqrt{2}}{2}(e_3-e_5)\ \ {\rm and}\ \
e_\pm(\gamma,\delta) = \sqrt{2}(\frac{\gamma}{2}(e_3+e_5) + \delta
e_4)\pm\sqrt{\frac{\gamma^2}{2}+\delta^2(e_3-e_5)}.$
Notice that the two values $e_\pm(\gamma,\delta)$ are equal when
$\frac{\gamma^2}{2}+\delta^2=0$. We have  
$f_{\gamma,\delta}(\zeta,\eta)=\zeta(\frac{\gamma}{2}f_5+\mu f_4)
+\eta f_3.$
Let us study the closures of these locally closed subsets.

If we project $B_X(\D)$ on  the grassmannian
$\mathbb{G}_Q(2,V_0)$ of totally isotropic subspaces of dimension $2$
of $V_0$, its image is the union of two lines meeting in a point
with the intersection point removed. 
The locally closed subset  $B_X(\D)$ is not connected and its closure
is made of two irreducible components $C_1$ and $C_2$. 

It is obvious that $B_X(\D''$) is closed and its image by the
projection onto $\mathbb{G}_Q(2,V_0)$ is made of two
points, thus $B_X(\D''$) has two connected components, isomorphic to
$\p^1$, $C_4$ and $C_5$. 

Finally, we can construct an isomorphism between $B_X(\D')$ and the
conic of isotropic points in $\p(\langle e_3,e_4,e_5\rangle)$ with two
points removed. Indeed, an element  $e(a,b,c)=a(e_3+e_5)+be_4+c(e_3-e_5)$
belongs to this conic if and only if  $2a^2+b^2+2c^2=0$. 
Set ${a=\frac{\sqrt{2}}{2}\gamma \textrm{ and }\
b={\sqrt{2}}\delta,}$ we have 
$c=\pm\sqrt{\frac{\gamma^2}{2}+\delta^2}.$
We obtain the following description
$$B_X(\D')=
\left\{\begin{array}{ll}
  \langle e_1\rangle
\subset
\langle e_1,e_2\rangle\subset \langle e_1,e_2,e(a,b,c)\rangle\\
\langle f_1\rangle\subset
\langle f_1, 2a f_2+b f_3\rangle \subset
\langle f_1,f_2,f_3\rangle
\end{array},
2a^2+b^2+2c^2=0,\ c\neq0
\right\}.\\ $$
The closure of $B_X(\D')$ is isomorphic to  $\p^1$ and gives us the component 
$C_3$. The intersections between the distinct components come from
previous descriptions. 

Let us go back to the general case. We know that there are $n$
admissible slicings.  
We denote them by  $\D^i$ for $i\in[1, n]$. If  $\D^1$ is the general one,
we know that the slicings $\D^i$ for $i>1$ are all slicings of the
diagram associated to the orbit 
 $\mathcal{O}_2$ in $\osp(2n-1,2n-2)$. The union of the corresponding
locally closed subsets is isomorphic to the fiber above the orbit
$\mathcal{O}_2$ in 
$\osp(2n-1,2n-2)$. To end the proof, we only have to check, by induction on $n$
that the locally closed subset $B_X(\D^1)$ which gives us the
components  $C_1$ and  $C_2$ meets the fiber as expected.  We proceed
as  in the $\osp(7,6)$ case.  
We begin by choosing a representative $u$ of the orbit $\mathcal{O}_2$
defined by $u(e_i)=f_{i-1}$ for 
$i\in[1,2n+1]\setminus\{n+2\}$ and $u(e_{n+2})=f_{n-1}$ (by convention $f_0=0$) 
where $(e_i)_{i\in[1,2n+1]}$ and $(f_i)_{i\in[1,2n]}$ are 
the bases of $V_0$ and $V_1$ respectively
such that $(e_i,e_j)=\delta_{i,2n+2-j}$ and 
$(f_i,f_j)=\delta_{i,2n+1-j}$. 
We can describe  the locally closed subsets $B_X(\D^i)$:
$$
\begin{array}{l}
  B_X(\D^1)=\left\{
    \begin{array}{ll}
\langle e_1\rangle
\subset
\langle e_1,\alpha e_2+\beta e_\pm\rangle \subset
\langle e_1,e_2,e_\pm\rangle\subset
\cdots\subset
\langle e_1,e_2,\cdots, e_{n-1},e_\pm\rangle \\
\langle \alpha f_1+\beta f_n\rangle\subset
\langle f_1,f_n\rangle\subset
\langle f_1,f_2,f_n\rangle \subset
\cdots\subset
\langle f_1,f_2,\cdots, f_n\rangle
\end{array},
\ \beta\neq 0
\right\}\\
B_X(\D^2)=\left\{
\begin{array}{ll}
  \langle e_1\rangle
\subset
\langle e_1,e_2\rangle\subset\langle e_1,e_2,\gamma e_3+\delta e_\pm\rangle\subset
\cdots
\langle e_1,e_2,\cdots,e_{n-1},e_\pm\rangle\\
\langle f_1\rangle\subset
\langle f_1, \gamma f_2+\delta f_3\rangle\subset
\langle f_1,f_2,f_3\rangle\subset
\cdots\subset
\langle f_1,f_2,\cdots, f_n\rangle
\end{array},\ \delta\neq0
\right\}\\ 
\end{array},$$
where 
$e_\pm=e_{n+1}\pm\frac{\sqrt{2}}{2}(e_n-e_{n+2}).$ The locally  closed subsets 
$B_X(\D^1)$ and $B_X(\D^2)$ have two connected components and their
closures have two irreducible components denoted   respectively by
$C_1$, $C_2$ for $B_X(\D^1)$ and  
$C_3$ $C_4$ for $B_X(\D^2)$. These components intersect each other as
predicted.

For the sets  $B_X(\D^i)$ with $i\geq2$, we only check that
the three first subspaces of the complete flag of $V_0$ are always 
$\langle e_1\rangle\subset 
\langle e_1e_2\rangle\subset
\langle e_1,e_2,e_3\rangle.$ In particular, the closures of the
locally closed  subsets 
$B_X(\D^i)$ for $i\geq2$ do not meet  the closure of $B_X(\D^1)$.
\end{preu}

\subsection{Orbit $\mathcal{O}_3$}

\begin{prop}
  \label{prop-3}
Let $n\geq 2$ and $X\in\osp(2n+1,2n)$ be in the orbit
  $\mathcal{O}_3$.
The reduced fiber  $\pi^{-1}(X)_{\textrm{r{e}d}}$ is the union
of $2n-2$ irreducible components 
$(C_i)_{i\in[1,n-1]}$ and $(S_i)_{i\in[1,n-1]}$ such that
\begin{itemize}
\item the components $(C_i)_{i\in[1,n-1]}$ are isomorphic to $\p^1$, 
\item the components $(S_i)_{i\in[1,n-1]}$ are isomorphic to the blow-up of $\p^2$ in two distinct points,  
\item the components $(C_{i})_{i\in[1,n-1]}$ form a chain,
\item the components $(S_{i})_{i\in[1,n-1]}$ form a chain,
two components intersect along a $\p^1$,
\item the two chains intersect each other in a point on $C_{n-1}\cap S_{n-1}$.
\end{itemize}
\end{prop}

\begin{preu}
We proceed by induction as in the previous proposition. We first study 
the $\osp(7,6)$ case. There are five admissible slicings 
$(\D^i)_{i\in[1,5]}$ of $D$. We use the bases
$(e_i)_{i\in[1,7]}$ and $(f_i)_{i\in[1,6]}$ defined in  the proof of 
the previous  proposition to choose
a representative of the orbit $u$,  whose matrix is the following one
$$\left(\begin{array}{ccccccc}
  0&1&0&0&0&0&0\\
  0&0&0&1&0&0&0\\
  0&0&0&0&0&0&0\\
  0&0&0&0&1&0&0\\
  0&0&0&0&0&1&0\\
  0&0&0&0&0&0&1\\
\end{array}\right).$$
The locally closed subsets
are:
$$\begin{array}{l}
  B_X(\D^1)=\left\{
    \begin{array}{ll}
\langle ae_1+be_3\rangle
\subset
\langle e_1,e_3\rangle \subset
\langle e_1,e_2,e_3\rangle\\
\langle \alpha (af_1+bf_3)+\beta f_4\rangle\subset
\langle f_1,\alpha bf_3+\beta f_4\rangle\subset\langle f_1,f_2,\alpha
bf_3+\beta f_4 \rangle\end{array},
\ b\neq 0
\right\}\\
B_X(\D^2)=\left\{
\begin{array}{ll}
  \langle e_1\rangle
\subset
\langle e_1,\alpha e_2+\beta e_5 \rangle\subset \langle e_1,e_2,e_5\rangle\\
\langle \alpha f_1+\beta f_4\rangle\subset
\langle f_1, f_4\rangle\subset
\langle f_1,f_2,f_4\rangle
\end{array},\ \beta\neq0
\right\}\\ 
\hskip 1.47 cm\coprod\left\{
\begin{array}{ll}
  \langle e_1\rangle
\subset
\langle e_1, e_3 \rangle\subset \langle e_1,e_2,e_3\rangle\\
\langle \alpha f_1+\beta f_4\rangle\subset
\langle f_1, f_4\rangle\subset
\langle f_1,f_2,f_4\rangle
\end{array},\ \beta\neq0
\right\}\\ 
B_X(\D^3)=\left\{
\begin{array}{ll}
  \langle e_1\rangle \subset 
\langle e_1,c e_2+d e_3 \rangle\subset
\langle e_1,e_2,e_3\rangle\\
\langle f_1\rangle\subset
\langle f_1, \gamma (cf_2+df_3)+\delta f_4\rangle\subset 
\langle f_1, f_2, \gamma d f_3+\delta f_4
\rangle
\end{array},\ d\neq0
\right\}\\ 
B_X(\D^4)=\left\{
\begin{array}{ll}
  \langle e_1\rangle
\subset
\langle e_1,e_2 \rangle\subset 
\langle e_1,e_2,\gamma^2 e_3-2\gamma\delta e_4-2\delta^2 e_5 \rangle\\
\langle f_1\rangle\subset
\langle f_1, \gamma f_2+\delta f_4\rangle\subset
\langle f_1,f_2,f_4\rangle
\end{array},\ \delta\neq0
\right\}\\ 
\hskip 1.47 cm\coprod\left\{
\begin{array}{ll}
  \langle e_1\rangle
\subset
\langle e_1,e_2 \rangle\subset 
\langle e_1,e_2,e_3 \rangle\\
\langle f_1\rangle\subset
\langle f_1, \gamma f_2+\delta f_4\rangle\subset
\langle f_1,f_2,f_4\rangle
\end{array},\ \delta\neq0
\right\}\\ 
B_X(\D^5)=\left\{
\begin{array}{ll}
  \langle e_1\rangle
\subset
\langle e_1,e_2 \rangle\subset 
\langle e_1,e_2,e_3 \rangle\\
\langle f_1\rangle\subset
\langle f_1, f_2\rangle\subset
\langle f_1,f_2,\zeta f_3+\eta f_4\rangle
\end{array},\ \delta\neq0
\right\}.\\ 
\end{array}$$
The closure of $B_X(\D^1)$  has in its boundary the two following lines
$$ \left\{
    \begin{array}{ll}
\langle e_1\rangle
\subset
\langle e_1,e_3\rangle \subset
\langle e_1,e_2,e_3\rangle\\
\langle \alpha f_1+\beta f_4\rangle\subset
\langle f_1, f_4\rangle\subset\langle f_1,f_2, f_4 \rangle\end{array}
\right\} 
\textrm{  and }
\left\{
\begin{array}{ll}
  \langle e_1\rangle
\subset
\langle e_1, e_3 \rangle\subset \langle e_1,e_2,e_3\rangle\\
\langle f_1\rangle\subset
\langle f_1,xf_3+y f_4\rangle\subset
\langle f_1,f_2,xf_3+yf_4\rangle
\end{array}
\right\}.\\ 
$$
The second one is an exceptional divisor and belongs to the closure
of $B_X(\D^3)$ when $c$ vanishes. The second exceptional divisor is obtained 
by putting $\beta=0$ in $B_X(\D^1)$. 
We can see then that the second component of $B_X(\D^2)$ is in the closure of
$B_X(\D^1)$. In the same way, the closure of $B_X(\D^3)$ has in its boundary 
two lines: 
$$ \left\{
    \begin{array}{ll}
\langle e_1\rangle
\subset
\langle e_1,e_3\rangle \subset
\langle e_1,e_2,e_3\rangle\\
\langle f_1\subset
\langle f_1,\gamma f_2+\delta f_4\rangle\subset
\langle f_1, f_2,f_4\rangle \rangle\end{array}
\right\} 
\textrm{  and  }
\left\{
\begin{array}{ll}
  \langle e_1\rangle
\subset
\langle e_1, e_3 \rangle\subset \langle e_1,e_2,e_3\rangle\\
\langle f_1\rangle\subset
\langle f_1,f_2\rangle\subset
\langle f_1,f_2,xf_3+yf_4\rangle
\end{array}
\right\}.\\ 
$$
The  second line is a component of $B_X(\D^5)$. The rest of the proposition 
follows from this description.  
For the general case, we proceed as for the orbit $\mathcal{O}_2$. We omit 
the details of the proof. 
\end{preu}

\end{document}